\documentclass{article}

\usepackage{arxiv}

\usepackage[utf8]{inputenc} 
\usepackage[T1]{fontenc}    
\usepackage{hyperref}       
\usepackage{url}            
\usepackage{booktabs}       
\usepackage{amsfonts}       
\usepackage{nicefrac}       
\usepackage{microtype}      
\usepackage{lipsum}		
\usepackage{graphicx}
\usepackage{natbib}
\usepackage{doi}
\usepackage{caption}
\usepackage{subcaption}
\usepackage{algorithm}
\usepackage{algorithmic}
\usepackage{amsmath}

\title{Multi-Axis and Multi-Vector Gradient Estimations: Using Multi-Sampled Complex Unit Vectors to Estimate Gradients of Real Functions}

\date{October 5, 2023}	

\author{ 
 \href{https://orcid.org/0000-0003-3618-4166}{\includegraphics[scale=0.06]{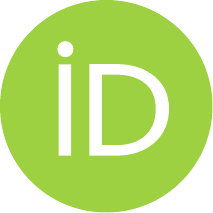}\hspace{1mm}Ergun Akleman}\thanks{Joint with Computer Science and Engineering Department.} \\
	Visual Computing \& Computational Media,\\ Texas A\&M University, College Station, TX, 77831\\
	\texttt{ergun@tamu.edu} \\
  \And
	 \href{https://orcid.org/0000-0002-3492-0628}{\includegraphics[scale=0.06]{orcid.pdf}\hspace{1mm}Alan D. Freed}\\
Department of Mechanical Engineering\\ 
Texas A\&M University, College Station, TX, 77831\\
	\texttt{afreed@tamu.edu } \\
}

\renewcommand{\shorttitle}{Multisample Gradient Estimation with Complex Vectors}

\hypersetup{
pdftitle={A template for the arxiv style},
pdfsubject={q-bio.NC, q-bio.QM},
pdfauthor={David S.~Hippocampus, Elias D.~Striatum},
pdfkeywords={First keyword, Second keyword, More},
}

\begin{document}
\maketitle

\begin{abstract}
In this preliminary study, we provide two methods for estimating the gradients of functions of real value. Both methods are built on derivative estimations that are calculated using the standard method or the Squire-Trapp method \cite{squire1998} for any given direction. Gradients are computed as the average of derivatives in
uniformly sampled directions. The first method uses a uniformly distributed set of axes that consists of orthogonal unit vectors that span the space. The second method only uses a uniformly distributed set of unit vectors. Both methods essentially minimize the error through an average of estimations to cancel error terms. Both methods are essentially a conceptual generalization of the method used to estimate normal fractal surfaces \cite{hart1989}. 
\end{abstract}  

\section{Introduction}

In this work, we study the gradients of general functions of the real value of $n$ dimensional $F:\Re^n \rightarrow \Re$. We will simply write these functions as
$F(\mathbf{p})$, where $\mathbf{p}$ is a point in a $n$-dimensional space, $\Re^n$. In other words, $\mathbf{p}=x$ in 1D, $\mathbf{p}=(x,y)$ in 2D, $\mathbf{p}=(x,y,z)$ in 3D, etc.  Our goal is to estimate the gradient of $F(\mathbf{p})$ numerically at any given point $\mathbf{p}_0$ in our $n$ dimensional space. We do not expect that the function has a well-defined derivative. In this short treatment of the problem, we focus only on the 2D and 3D gradients, but the methods are general and can be extended to higher dimensions. 

We essentially provide two types of estimation: (1) Multi-Axis and (2) Multi-Vector. Both are conceptual generalizations of the multi-sample method \cite{hart1989}, which was originally developed to compute normals on fractal surfaces. These are also closely related to the classical method. Note that the classical method is based on taking and combining the $n$ number of 1D derivatives in orthogonal directions for the gradients in the $n$ dimensional space. In multi-axis gradient estimation, we calculate the gradient in multiple orthonormal directions and use their average as an estimate. In multi-vector gradient estimation, we calculate the gradient in multiple directions and use their average as an estimate. 

Both of these methods require a good estimation of the derivatives in a given direction. To compute derivatives in the given direction, we consider using both the standard numerical approximation of derivatives, that is, finite differences, and Squire and Trapp's \cite{squire1998} method of using imaginary numbers.  These give us six possible ways to compute a gradient: (1) single-axis sampling using finite differences, (2) single-axis sampling using imaginary values, (3) multi-axis sampling using finite differences, (4) multi-axis sampling using imaginary values, (5) multi-vector sampling using finite differences, and (6) multi-vector sampling using imaginary values. Our initial qualitative assessment suggests that multi-axis using imaginary values provides the best results. Unlike the case of a single axis, in the multi-axis, many equally sampled orthogonal directions can provide a better estimate of the actual gradient as an average of multiple estimations of gradient vectors.  Multi-vector method also provides an estimate of the average gradient along a unit circle or a unit sphere (or a unit hypersphere for higher dimensions) around the given point \cite{hart1989}, however, since these directions may not include the orthogonal sets of vectors that represent randomly rotated axis, the error can be higher. 

\section{Preliminaries and Previous Work}

In this section, we will briefly discuss existing methods for estimating derivatives and gradients of functions whose derivatives may or may not exist. In these methods, it is sufficient that the function can work as a black box that gives us a single real number for any given point in the domain of the function. We are only interested in n-dimensional Euclidean domains that are defined as $\Re^n$, where $n$ is the dimension of the domain, which we call nD such as 1D, 2D, and 3D. In this study, we only discuss the estimation of gradients in 2D and 3D. Note that 1D gradients are simply derivatives. 

\subsection{Estimation of Derivatives}

The classical method for estimating the derivative of a 1D function $F(x)$ at a point $x_0$ is based on the standard central difference formula to define a derivative
$$\left. \frac{\mathrm{d} F(x)}{\mathrm{d} x} \right|_{x_0} = \lim_{h \rightarrow 0} \frac{F(x_0+h)-F(x_0-h)}{2h}$$
is as follows:
$$\left. \frac{\mathrm{d}F(x)}{\mathrm{d}x} \right|_{x_0} \approx \frac{F(x_0+h)-F(x_0-h)}{2h}$$

In theory, the error in this numerical calculation must be $O(h^2)$. However, the round-off error that occurs when subtracting the two terms prevents such accuracy in practice for small values of $h$. In 1998, Squire and Trapp \cite{squire1998} introduced a new way to improve the precision of estimation of derivatives of 1D functions with real value using complex algebra.  Squire and Trapp used a Taylor series to derive an alternative formula involving complex numbers and obtained the following approximation:
$$\left. \frac{\mathrm{d}F(x)}{\mathrm{d}x} \right|_{x_0} \approx \frac{ \Im ( F (x_0 + i h ))}{h} . $$
This method symbolically produces the same result as the standard derivative. To demonstrate this, consider that we have a general quadratic function $F(x) = a x^2+ bx + c$, we can ignore the constant term since it will not appear in the imaginary term, and therefore: 
\begin{eqnarray} 
\frac{dF}{dx} &=&\frac{1}{h} \Im \left( F(x+ih) \right)  \nonumber \\ 
&=& \Im \left( a (x+ih)^2 + b(x+ih) \right) \nonumber \\
&=& \frac{1}{h} \Im \left( a x^2+ 2 a x i h - a h + b x + i b h) \right) \nonumber \\ 
&=& \frac{1}{h} \left( 2 a x  h +  b h) \right) \nonumber \\ 
&=&  2 a x   +  b \nonumber \\ 
\end{eqnarray}
It is possible to demonstrate that for any function this method symbolically provides the same derivative as the original derivative. The advantage of this formula is that it not only gives higher precision digits for small values of $h$, but the precision stays the same for smaller values of $h$ \cite{greenwell1998}. 

Despite this advantage, this method is useful only for 1D cases. There is a need for a method that computes the gradient for cases of higher dimensions.  In this study, we generalize Squire and Trapp's method to handle higher-dimensional cases.

It is important to note that neither the standard method nor Squire and Trapp's imaginary method requires $F(x)$ to have a well-defined derivative. Even if $F(x)$ does not have a well-defined derivative at a given point, these methods provide an answer. This is very useful in Computer Graphics. 

\subsection{Multisample Method for Gradient Computation for Fractal Functions}
\label{sec_hart}

\begin{figure}
        \begin{subfigure}[t]{0.23\textwidth}
        \includegraphics[width=1.0\textwidth]{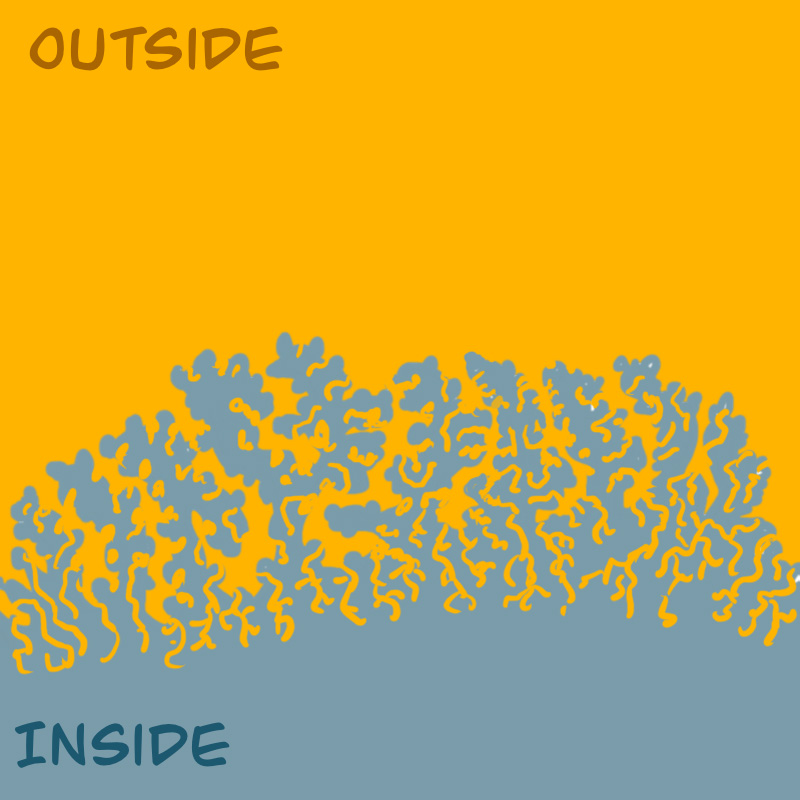}
        \caption{A black-box function that defines an interface. }
        \label{fig_MultiSample00}
    \end{subfigure}
    \hfill
        \begin{subfigure}[t]{0.23\textwidth}
        \includegraphics[width=1.0\textwidth]{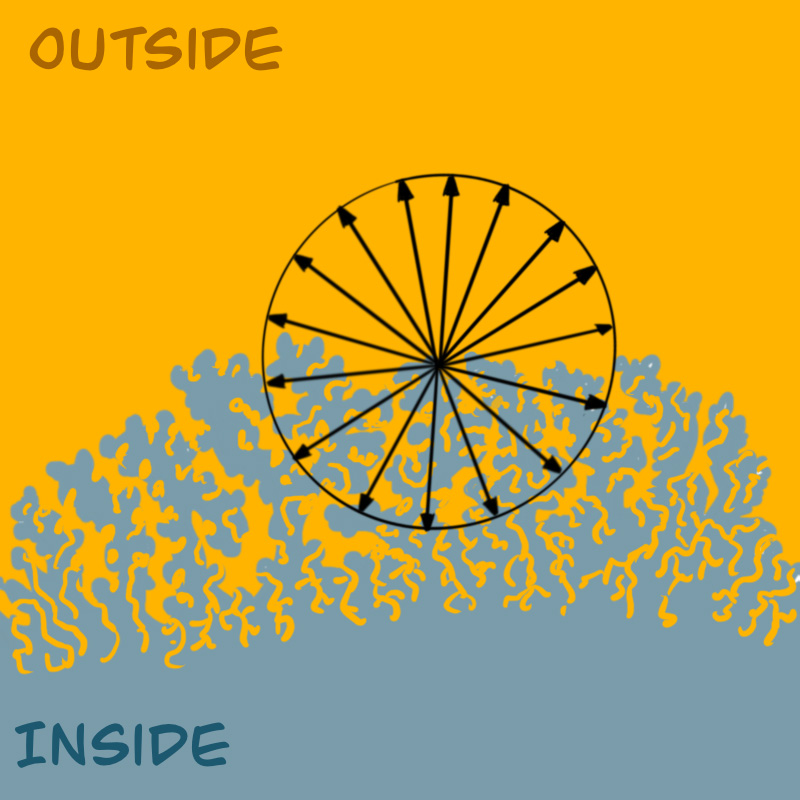}
        \caption{Uniformly distributed 2D unit vectors. }
        \label{fig_MultiSample01}
    \end{subfigure}
    \hfill
    \begin{subfigure}[t]{0.23\textwidth}
        \includegraphics[width=1.0\textwidth]{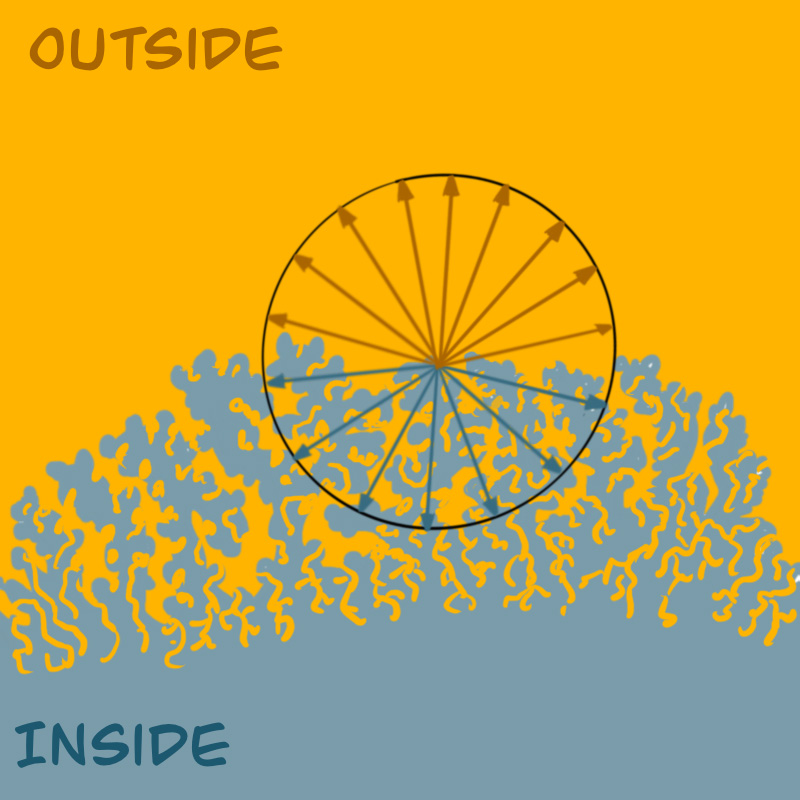}
        \caption{Classifications of unit vectors as inside and outside. }
        \label{fig_MultiSample02}
    \end{subfigure}
    \hfill
        \begin{subfigure}[t]{0.23\textwidth}
        \includegraphics[width=1.0\textwidth]{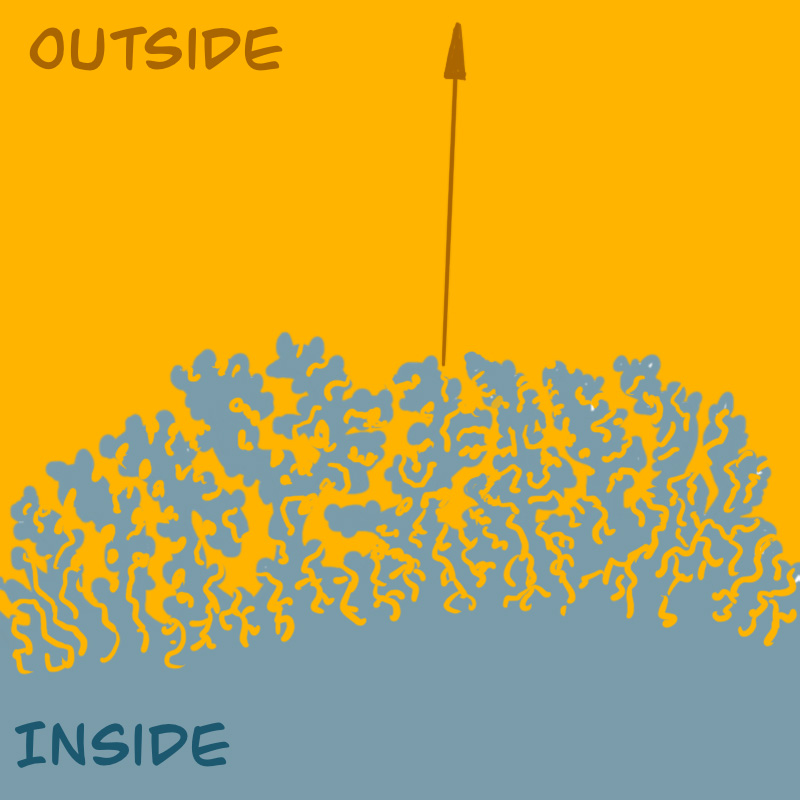}
        \caption{An estimated normal by adding outside vectors. }
        \label{fig_MultiSample03}
    \end{subfigure}
    \hfill
        \caption{Text.   }
        \label{fig_MultiSample0}
\end{figure}

Hart et al. introduced the multi-sample method to render 3D fractals. To compute diffuse illumination and specular reflection, one needs to estimate the surface normal in regions around the boundaries of 3D fractals \cite{hart1989}. The fractals are not represented by nice analytical functions, and the boundaries have a very high curvature; the estimation is even harder. We observe that Hart's solution is to view the fractal description algorithms in a black-box function that can provide an inside-outside test. Figure~\ref{fig_MultiSample00} shows an example of the boundaries defined by such black-box functions as an interface between inside and outside. Such a black-box function $F(\mathbf{p})$ can be given in the following form: 
$$ 
F(\mathbf{p}) =\left\{  
\begin{array}{ll} 
      a &
      \mbox{if $\mathbf{p}$ is outside} \\ 
      b &
      \mbox{if $\mathbf{p}$ is inside} \\ 
\end{array} 
\right.    	            
$$ 

where $a$ and $b$ are real numbers with $a>b$ to be used to identify inside and outside. 

Note that this is a terrible case. The function is not only discontinuous, but its boundary is so complicated that it is hard to assign a normal vector. To compute the gradient in such cases, the first step in Hart's method is to create a set of unit vectors $\vec{\mathbf{n}}_k$ that are uniformly distributed on the surface of a unit sphere of dimensions $n$, where $k=0,1,\cdots, K-1$ (see Figure~\ref{fig_MultiSample01}). We then classify these vectors as inside- and outside-vectors for a given position $\mathbf{p}_0$. The vector $\vec{\mathbf{n}}_k$ is an inside vector if $F(\mathbf{p}_0 + \vec{\mathbf{n}}_k) = a$, otherwise it is an outside vector (see Figure~\ref{fig_MultiSample02}). The normal vector is then calculated by adding all inside vectors (see Figure~\ref{fig_MultiSample03}). Note that the main idea in this method is that the tangent components cancel out each other, while the perpendicular directions reinforce each other. If we choose $a=1$ and $b=0$, this can be formulated as follows. 
$$
\left. \nabla F(\mathbf{p}) \right|_{\mathbf{p}=\mathbf{p}_0}  \approx \sum_{k=0}^{K-1}  F(\mathbf{p}_0 + \vec{\mathbf{n}}_k) \vec{\mathbf{n}}_k. \nonumber
$$
Note that the values of $a$ and $b$ can be anything (as far as the $a>b$ condition is satisfied) to obtain a somewhat acceptable normal vector. Since this uses all vectors, the result is slightly more robust. A good choice for the values of $a$ and $b$ is $a=1$ and $b=-1$. 

\subsection{Antialising with Multisample Method} 

One problem with the multisampling method is that it can suffer from aliasing. To solve it, we need to randomize the unit vectors while keeping their distribution intact \cite{cook1984distributed,cook1986stochastic}. The simple solution is to rotate all vectors at random. For each new position to compute the gradient vector, we can create a random rotation matrix $R$ in $SO(3)$, then update all unit vectors as follows: 
\begin{equation} 
\vec{\mathbf{n}}_k  \leftarrow R \vec{\mathbf{n}}_k. \nonumber
\end{equation}
Note that in this case, we cannot randomize vectors since that can change the distribution. The random rotation of all vectors creates a new randomization while keeping the original distribution of vectors. This idea is nothing more than a reformulation of the original stochastic sampling \cite{cook1984distributed,cook1986stochastic} that is adapted for a uniformly distributed set of unit vectors. 

\section{Theoretical Framework}

Now, we are ready to extend the multi-sample method to be used beyond to estimate the gradients of Fractal-type black-box functions. First, we need to define the derivative in a given direction. 

\subsection{Derivative in a Given Direction}

Numerical derivatives in 1D can be easily extended to compute the derivative of any high-dimensional function in any given direction. Let $\vec{\mathbf{n}}$ denote a unit vector, and $\mathbf{p}_0$ denote a position in a $n$-dimensional space. The $F(\mathbf{p})$ becomes a one-dimensional function along the direction of a ray $\vec{n}$ emanating from $\mathbf{p}_0$.  The equation for this ray can be written simply using an auxiliary parameter $t$ as
$$\mathbf{p} = \mathbf{p}_0 + \vec{\mathbf{n}} t $$
and a one-dimensional function along this ray is given in terms of the parameter $t$ as
$$ F(\mathbf{p}_0, \vec{\mathbf{n}}(t) )= F(\mathbf{p}_0 + \vec{\mathbf{n}} t)$$
whose derivative along this ray can be simply estimated by the classical method as follows
\begin{equation}
\left. \frac{dF(\mathbf{p} + \vec{\mathbf{n}}t)}{dt} \right|_{\mathbf{p}=\mathbf{p}_0} \approx \frac{F(\mathbf{p}_0 + \vec{\mathbf{n}}h)-F(\mathbf{p}_0 - \vec{\mathbf{n}}h)}{2h} \label{Eq:Classical}
\end{equation}
We can also use Squire and Trapp's method by replacing $t$ with $i h$ \cite{almohy2010,squire1998} as follows
\begin{equation}\left. \frac{dF(\mathbf{p} + \vec{\mathbf{n}}t)}{dt} \right|_{\mathbf{p}=\mathbf{p}_0} \approx \frac{ \Im ( F (\mathbf{p}_0 + \vec{\mathbf{n}}_k i h ))}{h} .  \label{Eq:SquireTrapp}
\end{equation}

\subsection{Derivative in a Given Line} 

Note that a line passing through a point consists of two rays defined by two opposing vectors. The derivative has to be the same in both directions once we consider the vector directions. In other words, we expect the following.  
$$ F(\mathbf{p}_0, \vec{\mathbf{n}}(t)) \vec{\mathbf{n}} = F(\mathbf{p}_0,\vec{\mathbf{-n}}(t) ) \vec{\mathbf{-n}} $$
where $\vec{\mathbf{n}} + (\vec{\mathbf{-n})} = 0$. This is a useful property, since if we estimate the derivative in both directions, their average can provide us with a better estimate by averaging errors if we assume the errors to be random. Note that Equation~\ref{Eq:Classical} is written as an average of the two derivatives. However, the negative creates a problem in estimation. Equation~\ref{Eq:SquireTrapp} does not have this problem, and the result can be improved by computing the average as follows: 

\begin{equation}\left. \frac{dF(\mathbf{p} + \vec{\mathbf{n}}t)}{dt} \right|_{\mathbf{p}=\mathbf{p}_0} \approx \frac{ \Im ( F (\mathbf{p}_0 + \vec{\mathbf{n}}_k i h )) + \Im ( F (\mathbf{p}_0 + (\vec{\mathbf{-n}}_k i h )))}{2h} .  \label{Eq:SquireTrapp2}
\end{equation}

\section{Multi-Axis Gradient Estimation} 

Now, we are ready to develop the first method, which we call the multi-axis gradient. In this case, we simply compute the gradient for multiple axes and find the average. The following subsection presents the multi-axis gradient method by estimating the gradient in multiple sets of orthonormal vectors. 

\subsection{Gradient with Orthonormal Vectors} 

Let $\vec{\mathbf{n}}_j$, where $j=0,1,\cdots, n-1$  denote any set of orthonormal vectors in the $n$ dimensional space, $\Re^n$. The orthonormal means that all these unit vectors are perpendicular to each other.  

$$ 
\nabla F(\mathbf{p})|_{\mathbf{p}=\mathbf{p}_0} = \sum_{j=0}^{n-1} \vec{\mathbf{n}}_j \left. \frac{dF(\mathbf{p} + \vec{\mathbf{n}_j}t)}{dt} \right|_{\mathbf{p}=\mathbf{p}_0}. 
$$
Our intuition observation is that this equation must be correct regardless of how orthogonal vectors are chosen. Now, consider a 2D case to test this intuitive observation. 
$$ 
\nabla F(\mathbf{p}) \approx \vec{\mathbf{n}_0}  \frac{dF(\mathbf{p} + \vec{\mathbf{n}_0}t)}{dt}  + 
\vec{\mathbf{n}_1}  \frac{dF(\mathbf{p} + \vec{\mathbf{n}_1}t)}{dt} . 
$$
where two perpendicular unit vectors can be given either $\vec{\mathbf{n}}_0 = (\cos \theta, \sin \theta )$ and $\vec{\mathbf{n}}_1 = ( -\sin \theta, \cos \theta )$, or $\vec{\mathbf{n}}_0 = (\cos \theta, \sin \theta )$ and $\vec{\mathbf{n}}_1 = ( \sin \theta, -\cos \theta )$. 
Now, let us check if this gradient gives us the actual formula of the gradient: 
\begin{eqnarray} 
\nabla F(\mathbf{p})  & = &  \frac{dF(\mathbf{p} + \vec{\mathbf{n_0}}t)}{dt} \vec{\mathbf{n_0}}+
 \frac{dF(\mathbf{p} + \vec{\mathbf{n_1}}t)}{dt} \vec{\mathbf{n_1}} \nonumber\\ 
& = &  \frac{dF(x+t\cos \theta, y+ t\sin \theta )}{dt} (\cos \theta, \sin \theta )+
 \frac{dF(x-t\sin \theta, y+ t\cos \theta )}{dt} (-\sin \theta,\cos \theta ) \nonumber \\ 
& = &  \left(\frac{\partial F(x,y)}{\partial x} \cos \theta + \frac{\partial F(x,y)}{\partial y} \sin \theta \right) (\cos \theta, \sin \theta ) +
\left( \frac{-\partial F(x,y)}{\partial x} \sin \theta + \frac{\partial F(x,y)}{\partial y} \cos \theta \right) (-\sin \theta, \cos \theta ) \nonumber \\
&= & \frac{\partial F(x,y)}{\partial x} (1,0) +\frac{\partial F(x,y)}{\partial y}  (0,1)\nonumber \\
\end{eqnarray}
Note that this is not complete proof since it just holds for 2D cases. On the other hand, it provides additional support to our intuitive observation. 

Let us assume that our intuitive observation supported for 2D cases is correct for higher dimensions. This means that we can compute multiple estimations of the gradient using different sets of orthonormal vectors. For calculation, we can use either the classical method or the imaginary method such as the one given below: 
$$ \nabla F(\mathbf{p})|_{\mathbf{p}=\mathbf{p}_0} \approx \frac{1}{h} \sum_{j=0}^{n-1} \Im \bigl( F (\mathbf{p}_0 + \vec{\mathbf{n}}_j i h ) \bigr) \vec{\mathbf{n}}_j . $$
Each of these calculations provides an estimate based on the precision of the computation method. Now, let us also assume that in these estimates the errors are symmetric. we can calculate the gradient as the average of multiple estimates. If the symmetry condition holds, this average should give us a better estimate of the gradient. An additional property is that, even when derivatives do not exist, it gives us a good estimate of the gradient vector. 

It is important to note that any orthonormal set of vectors in n-dimension can be obtained using a transformation matrix $M \in O(n)$. Here, $O(n)$ is called the orthogonal group in n-dimension, and it is the group of all distance-preserving transformations of an n-dimensional Euclidean that preserve a fixed point. Note that these operations include mirror operations along with rotation operations. Now, let $M_k \in O(n)$ where $k=0,1,\cdots, K-1$ and let $\vec{\mathbf{n}}_j$ denote a set of initial orthonormal vectors where $j=0,1,\cdots, n-1$, then we can compute the average gradient starting from the following formula:
$$ 
\nabla F(\mathbf{p})|_{\mathbf{p}=\mathbf{p}_0} = \frac{1}{K} \sum_{k=0}^{K-1} \sum_{j=0}^{n-1} M_k\vec{\mathbf{n}}_j \left. \frac{dF(\mathbf{p} + M_k\vec{\mathbf{n}_j}t)}{dt} \right|_{\mathbf{p}=\mathbf{p}_0}. 
$$

To ensure that the errors are symmetric, we need to get $M$'s by equally sampling $O(n)$. If the samples are skewed in one region, the average may not necessarily be good. Higher dimensions require special care to obtain the symmetric error. 

We want to point out that the term
$$\left. \frac{dF(\mathbf{p} + M_k\vec{\mathbf{n}_j}t)}{dt} \right|_{\mathbf{p}=\mathbf{p}_0}$$
can be calculated by either the standard method or the imaginary method to obtain an estimate. This gives us two ways to compute the gradient. 

\subsection{Demonstration of Multi-Axis Method for 2D}
\label{sec_2DMulti-Axis}

It is easy to demonstrate this work for 2D, as it is straightforward to sample equally $O(2)$. Consider a set of unit vectors that are equally distributed in a circle. These vectors can be created using the functions $\cos$ and $\sin$ as follows:
$$\vec{\mathbf{n}}_k = \left( \cos \frac{\pi k }{K}, \sin \frac{\pi k }{K}  \right) , 
\qquad k = 0 , 1 , 2 , \ldots , K-1.$$
Here, $M$'s are implicit, which are given as rotation matrices that are equally spaced $S(2)$. 
If $K$ is divisible by four, then this set of vectors provides us with the $K/2$ number of the axis. 
For any given point $\mathbf{p}=(x,y)$ in any given direction $\vec{\mathbf{n}}_k$, we approximate the derivative along this ray as:
\begin{eqnarray}
x(t) &=& x + t \cos \frac{\pi k }{K} \nonumber \\
y(t) &=& y + t\sin \frac{\pi k }{K} \nonumber 
\end{eqnarray}
so that for any given 2D function $F(x,y)$, the derivative along this ray is computed as follows by replacing $t$ with $ih$: 
$$\nabla F(\mathbf{p})  \approx \frac{ \Im( F (x_0 + \cos \frac{ \pi k }{K} i h  ,  y_0 + \sin \frac{ \pi k}{K}  i h ))}{h} .$$

To see that this method symbolically produces the same result as the standard gradient, now consider the simple circle equation $F(x,y) = x^2+y^2$, we can ignore the constant term since it will not appear in the imaginary term, and therefore: 
\begin{equation} 
\frac{1}{h} \Im \left( \left(x + \cos \frac{ \pi k }{K}  i h \right)^2 + \left( y + \sin \frac{ \pi k}{K}  i h \right)^2 \right)  =
 2 x \cos \frac{ \pi k }{K}  + 2 y \sin \frac{ \pi k}{K} 
\end{equation}
whose gradient is computed as: 
\begin{equation} 
\nabla F(\mathbf{p})  \approx \sum_{k=0}^{K-1} \left( 2 x \cos \frac{\pi k }{K}  + 2 y \sin \frac{ \pi k}{K} \right)  \left( \cos \frac{ \pi k }{K}, \sin \frac{ \pi k }{K}  \right) .
\end{equation}

For $K=4$, the formula consists of the addition of the four vectors as follows: 
\begin{eqnarray} 
\nabla F(\mathbf{p}) &\approx& \sum_{k=0}^{3} \left( 2 x \cos \frac{\pi k }{4}  + 2 y \sin \frac{ \pi k}{4} \right)  \left( \cos \frac{ \pi k }{4}, \sin \frac{ \pi k }{2}  \right) \nonumber \\
&=&  \left( 4 x \right)  \left( 1, 0  \right)  +   \left( 4 y \right)  \left( 0, 1  \right) \nonumber \\
&=& (4x, 4y)
\end{eqnarray}
In this case, there are two sets of axes, and when we divide this result by the number of axes, the computation gives us an actual gradient $(2x, 2y)$. This is the simplest version of the multi-axis. For any value of $K$ divisible by $4$, we obtain a multi-axis method. In this case, the result becomes 
\begin{equation} 
\nabla F(\mathbf{p}) \approx \left( 2 x  \sum_{k=0}^{K-1}  \cos^2 \frac{\pi k }{K}  + 2 y \sum_{k=0}^{K-1} \sin \frac{ \pi k}{K} \cos \frac{ \pi k }{K} , 
2 y \sum_{k=0}^{K-1}   \sin^2 \frac{\pi k }{K}  +  2 x \sum_{k=0}^{K-1} \sin \frac{ \pi k}{K} \cos \frac{ \pi k }{K}\right)  \label{Equation_1}
\end{equation}
In Equation \ref{Equation_1}, if $K$ is even 
\begin{eqnarray} 
\sum_{k=0}^{K-1}  \cos^2 \frac{\pi k }{K}  &=& \frac{K}{2}  \label{Equation_2} \\
\sum_{k=0}^{K-1}  \sin \frac{ \pi k}{K} \cos \frac{ \pi k }{K} &=& 0  \label{Equation_3}
\end{eqnarray}
As a result, the Equation \ref{Equation_1} is calculated as 
\begin{equation} 
\nabla F(\mathbf{p}) \approx (Kx, Ky)  \label{Equation_4}
\end{equation}
Note that to obtain the actual gradient, we need to take the average by dividing the resulting vector with $K/2$, which corresponds to averaging of the computation. Therefore, the multi-axis method can potentially give better numerical results since the average removes the errors. 

We also point out that in this example $K$ must be divisible by four since we use uniformly spaced vectors. For single-axis computation, we will need to use one of the two following axes:
\begin{eqnarray}
\vec{\mathbf{n}}_0 = \left(\cos \theta, \sin \theta \right) \; \; \; \vec{\mathbf{n}}_1 =\left( \cos (\theta+\pi/2), \sin (\theta+\pi/2) \right) \nonumber \\
\vec{\mathbf{n}}_0 =\left( \cos \theta, \sin \theta \right) \; \; \; \vec{\mathbf{n}}_1 =\ \left( \cos (\theta-\pi/2), \sin (\theta-\pi/2) \right)
\end{eqnarray} 
For higher dimensions, $K$ must be at least a multiple of $n$ to have a multi-axis. In higher dimensions, the normalization term will be $K/n$ instead of $K/2$. It is also important to note that this symbolic computation is possible only for the $2D$ case. For higher dimensions, we do not always get such nice closed formulas. Another important point is that even in the 2D case, there will be round-off errors that can accumulate in numerical computations. 

Another tricky problem is the number of samples. Note that to estimate gradients with the single-axis method, we only need the $n$ number of orthonormal vectors in $n$ dimensions. The standard finite-difference formula effectively uses its mirror. In practice, we use the $2N$ number of vectors. In other words, the axis consists of $2n$ vectors.  This is the reason in a 2D space that we need $4$ vectors instead of $2$. If we populate the space with these axes, we need to have $K=2n M$ number of vectors, where any $M$ is any positive integer. This is why we need $K$ divisible by $4$ in 2D. In higher dimensions, populating the space with a uniformly distributed set of $M$ number of axes may not always be straightforward. 

Now, assume that $K$ is not divisible by $2n$ or even $n$. We also assume that the vectors may not include orthonormal sets. It is still possible to calculate a sum and divide the sum by $K/n$ to obtain an average that gives us an estimation of the gradient. The main problem is that we will not have the theoretical guarantee that there will be a $K/n$ number of axes, and each of them gives us a good estimate of the actual gradient. On the other hand, from the multi-sample method discussed in Section~\ref{sec_hart} we know that it could still provide an approximate estimate. This observation gives us the next method. 

\subsection{Multi-Vector Gradient Estimation}

Now, we create a set of unit vectors $\vec{\mathbf{n}}_k$ that are uniformly distributed on the surface of a unit sphere of dimensions $n$, where $k=0,1,\cdots, K-1$.  This is an interesting problem in itself.  We assume that we can do this in the general case.  This problem is very simple in 2D as we discussed in the previous section. The problem is still reasonably easy in 3D regular and semi-regular polyhedra.  However, for higher-dimensional cases, we need to use regular and semi-regular polytops \cite{coxeter1973regular,coxeter1985regular}, which can be more difficult to conceptualize. Now, let us assume that it is possible to obtain uniform distribution in any dimension. 

For every unit vector $\vec{\mathbf{n}}_k$, we compute a derivative of $F(\mathbf{p})$ along the direction of a ray $\vec{\mathbf{n}}_k$ emanating from $\mathbf{p}_0$.  The equation for this ray can be written simply using an auxiliary parameter $t$ as 
$$\mathbf{p} = \mathbf{p}_0 + \vec{\mathbf{n}}_k t $$
whose derivative along this ray we estimate by replacing $t$ with $i h$ \cite{almohy2010,squire1998}  
$$\frac{ \Im ( F (\mathbf{p}_0 + \vec{\mathbf{n}}_k i h ))}{h} . $$
The gradient of $F(\mathbf{p})$ is then approximated as the sum of all derivative estimates multiplied with their corresponding unit vectors as follows: 
$$ \nabla F(\mathbf{p})|_{\mathbf{p}=\mathbf{p}_0} = \frac{1}{h} \sum_{k=0}^{K-1} \Im \bigl( F (\mathbf{p}_0 + \vec{\mathbf{n}}_k i h ) \bigr) \vec{\mathbf{n}}_k . $$
Note that this formula needs to be resized by dividing $K/n$ to compute the average over all the sample points. 

\section{Uniformly Distributed Vector Creation}

In the 2D case, the creation of equally spaced unit vectors can be done using the vertex positions of the regular polygons that are circumscribed inside a unit circle whose center is located in the origin such as a regular triangle, square, or regular pentagon. As discussed in Section~\ref{sec_2DMulti-Axis}, for the multi-axis method we need regular polygons with $4n$ sides such as square, regular octagon, regular dodecagon, and regular hexadecagon. Note that for the multi-axis method, it is also possible to use only half of the circle by eliminating the inverse vectors. In other words, if $\vec{\mathbf{n}}_k$ is used, $-\vec{\mathbf{n}}_k$ will not be used. This removes half of the vectors, and using only even numbers of vectors, we can estimate the gradient with the multi-axis method. 
Note that polygons with an even number of corners but not divisible by $4$, such as regular hexagons and regular decagons, cannot still be used for the multi-axis method since they do not include perpendicular vectors. Of course, these and regular polygons with odd numbers of sides can be used with the multi-vector method. This discussion provides a glimpse of how the problem can be complicated in higher dimensions.  

\begin{figure}
        \begin{subfigure}[t]{0.23\textwidth}
        \includegraphics[width=1.0\textwidth]{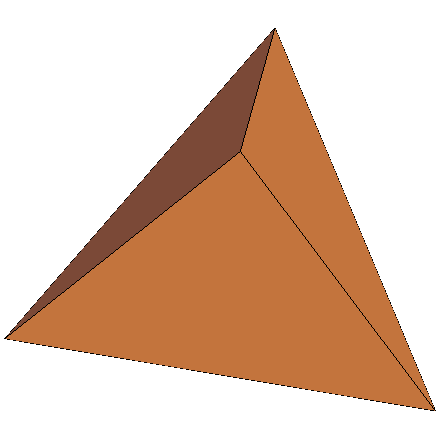}
        \caption{Tetrahedron. }
        \label{fig_tetrahedron.png}
    \end{subfigure}
    \hfill
            \begin{subfigure}[t]{0.23\textwidth}
        \includegraphics[width=1.0\textwidth]{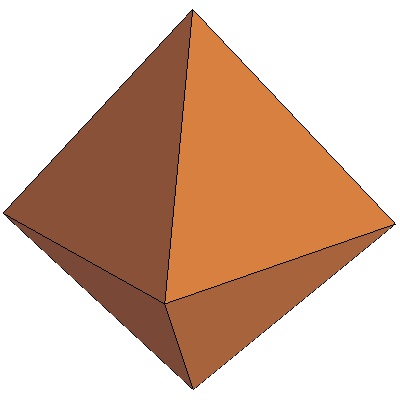}
        \caption{Octahedron. }
        \label{fig_octahedron.png}
    \end{subfigure}
    \hfill
    \begin{subfigure}[t]{0.23\textwidth}
        \includegraphics[width=1.0\textwidth]{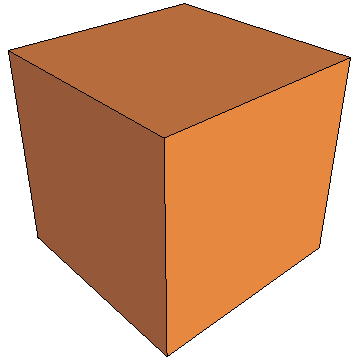}
        \caption{Cube. }
        \label{fig_cube.png}
    \end{subfigure}
    \hfill
        \begin{subfigure}[t]{0.23\textwidth}
        \includegraphics[width=1.0\textwidth]{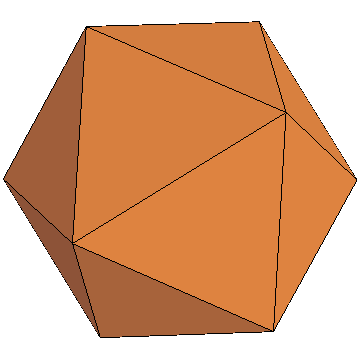}
        \caption{Icosahedron. }
        \label{fig_icosahedron.png}
    \end{subfigure}
    \hfill
            \begin{subfigure}[t]{0.23\textwidth}
        \includegraphics[width=1.0\textwidth]{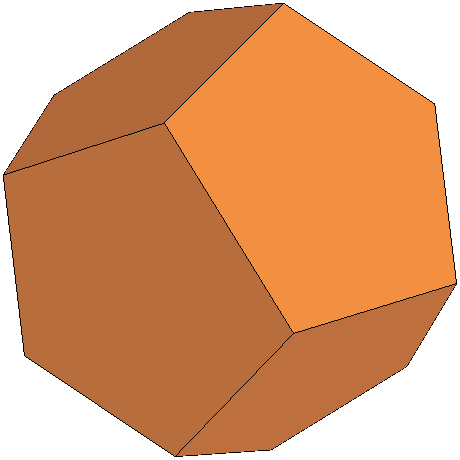}
        \caption{Dodecahedron. }
        \label{fig_dodecahedron.png}
    \end{subfigure}
    \hfill
    \begin{subfigure}[t]{0.23\textwidth}
        \includegraphics[width=1.0\textwidth]{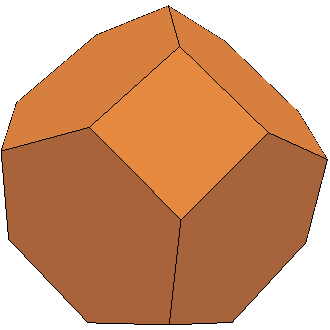}
        \caption{Truncated Octahedron. }
        \label{fig_truncatedoctahedron.png}
    \end{subfigure}
    \hfill
    \begin{subfigure}[t]{0.23\textwidth}
        \includegraphics[width=1.0\textwidth]{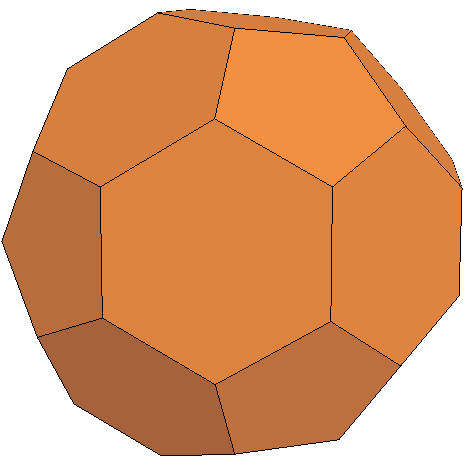}
        \caption{Soccer Ball. }
        \label{fig_soccerball.png}
    \end{subfigure}
    \hfill
    \begin{subfigure}[t]{0.23\textwidth}
        \includegraphics[width=1.0\textwidth]{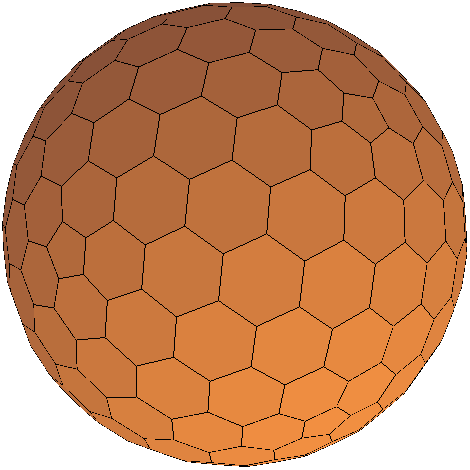}
        \caption{Geodesic Dome (Not Semi-Regular). }
        \label{fig_dome.png}
    \end{subfigure}
    \hfill
        \caption{Examples of regular and semiregular polyhedra. To use them for multisampling, their centers must be the origin, and every vertex position must be given by a unit vector. That unit vector provides the direction. 
        Note that not all of them can still be used for 3D multisampling. In these cases, the octahedron is the shape that corresponds to the classical method by providing three orthogonal directions. Note that the four vertices of the tetrahedron provide the same distribution as the vertices of the cube. Therefore, the tetrahedron is practically useless. The icosahedron provides 10 direction vectors, and the dodecahedron provides 6 direction vectors. For semiregular polyhedra, all vertices must be congruent for multisampling. The two examples here, the truncated octahedron and the soccer ball, are such examples. If we remove the same distribution criteria, we can use a wide variety of shapes such as the Geodesic Domes.   }
        \label{fig_regularpolyhedron}
\end{figure}

In the 3D case, the creation of equally spaced unit vectors can be done using the vertex positions of the regular or semiregular polyhedra that are circumscribed inside of a unit sphere whose center is located in origin. Some examples are shown in Figure~\ref{fig_regularpolyhedron}. However, most of them are useful only for multi-vector gradient estimation. For example, the cube cannot be used for multi-axis gradient estimation since its number of vertices is not divisible by $6$. 
The octahedron is the simplest that can be used for multi-axis gradient estimation (see Figure~\ref{fig_octahedron.png}). The icosahedron with its $12$ vertices can also be used for multi-axis gradient estimation (see Figure~\ref{fig_icosahedron.png}). Note that the number of vertices of the icosahedron is not only divisible by $6$, but its vertex positions consist of two distinct sets of perpendicular directions. Similarly, the truncated octahedron shown in Figure~\ref{fig_truncatedoctahedron.png} consists of four distinct sets of perpendicular directions and can be used for multi-axis gradient estimation. Any of these polyhedra. i.e. the octahedron, the icosahedron, and the truncated octahedron can be used to create new vectors by randomly rotating them. Since the initial distribution is uniform, the random rotations create vectors that can cover the surface of the sphere in a reasonably uniform way. We can approach this problem in higher dimensions in a similar way. 

\subsection{Discussion}

In this section, we provide a brief discussion on difficulty in comparison of the numerical computational methods with ground truths.  Let the shapes $S$ be defined through implicit functions $F(\mathbf{p})$'s as 
$$S = \{\mathbf{p} | F( \mathbf{p} ) \leq 0 \} . $$
One tricky problem for implicit functions, there exist infinitely many functions that can describe the same shape. Concerning gradients, we are interested in accurately describing tangent hyper-planes on our shape boundaries whose unit vectors are 
$$\vec{\mathbf{n}} = \frac{\nabla F( \mathbf{p} )}{\|\nabla F(\mathbf{p})\|}.$$
In other words, we are interested only in the direction of these normal vectors, independent of how the functions are formed. A simple example that demonstrates this idea is a circle. 

Notice that a circle of radius $r$ may be described using any of the following functions: 
\begin{eqnarray}
F_0(\mathbf{p}) &=& x^2 + y^2 - r^2 \nonumber \\
F_1(\mathbf{p}) &=& \frac{x}{r}^2 +  \frac{y}{r}^2 - 1 \nonumber \\
F_2(\mathbf{p}) &=& \sqrt{x^2 + y^2} - r \nonumber \\
F_3(\mathbf{p}) &=& (x^2 + y^2)^a- r^2a \nonumber 
\end{eqnarray}
where $a$ is an arbitrary positive real number. The lengths of the gradient vectors for these functions are vastly different, as can be easily seen in the following list: 
\begin{eqnarray}
\nabla F_0(\mathbf{p}) &=& 2 (x,y) \nonumber \\
\nabla F_1(\mathbf{p}) &=& \frac{2 (x, y) }{r^2} \nonumber \\
\nabla F_2(\mathbf{p}) &=& \frac{(x , y)}{\sqrt{x^2 + y^2}} \nonumber \\
\nabla F_3(\mathbf{p}) &=& 2a (x , y ) (x^2 + y^2)^{a-1} \nonumber
\end{eqnarray}
On the other hand, it is also easy to check that their unit gradients are the same on the boundary; they being precisely $\nabla F_2(\mathbf{p})$. 

Therefore, it is possible to claim that it is important to check the accuracy against the unit gradient instead of against its length. Note that if we agree with this argument, there is no need to divide $K/n$. We can simply normalize the gradient. One problem is that this normalization will also include some errors. Let us assume that we can ignore this and we want to compare numerical results with ground truth. In 2D, it can be easy  to visualize the precision of the unit numerical gradient $\vec{\mathbf{n}}_n$ with that of the ground truth $\vec{\mathbf{n}}_{gt}$ using $\cos \theta$, where $\theta$ is simply the angle between the two vectors, which can be computed using the dot product as follows:
$$\cos \theta = \vec{\mathbf{n}}_n \cdot \vec{\mathbf{n}}_{gt} . $$

A nice property of this error measure is that it is bounded between $[-1,1]$, where $1$ means no error, and $-1$ means the maximum possible. This error can be visually represented through gray-level colors by mapping $[-1,1]$ to $[0,1]$ using the following transformation:
$$\epsilon = \frac{\cos \theta + 1}{2}$$

Based upon these observations, we initially conjecture that it can be possible to compare the numerical gradient computations with the classical numerical gradient computations and with the ground truth. One simply creates an ancillary image of some image whose gradients are sought.  The points comprising this ancillary image are dot products that describe an error, thereby producing a visual representation of the error over its spatial domain.  Images can be constructed using classical methods to estimate numerical gradients. If a method turns out to be a significantly brighter image, then we can visually claim that that method is superior for the problem considered.

Unfortunately, in these cases, we quickly realized that even ground truths are approximate since they should be calculated from symbolic results. The additional normalization term makes the ground-truth computation not that reliable. We also realized that low-dynamic range images are not suitable for comparison. All methods provide relatively bright results that can be difficult to visually compare results in a large domain. Therefore, we provide a qualitative assessment of these methods. 

\begin{enumerate}
\item \textbf{Single-Axis/Finite-Difference:} This approach potentially suffers from inaccuracies that come from a low number of sampling directions due to the use of only a single set of orthogonal directions, and from the roundoff error caused by the difference in the finite-difference formul\ae that is used to estimate the derivative.    
\item \textbf{Single-Axis/Imaginary:}   The imaginary step derivative no longer succumbs to the roundoff error, but the results still suffer from inaccuracies due to the low number of sampling directions, because only a single orthogonal direction is used.
\item \textbf{Multi-Axis/Finite-Difference:} The integral method eliminates the problem of low numbers of sampling directions, but the roundoff error from the standard derivative approximation can still be significant.    
\item \textbf{Multi-Axis/Imaginary:}  The imaginary step derivative coupled with integral sampling significantly improves the results.
\item \textbf{Single-Axis/Imaginary:}   The imaginary step derivative no longer succumbs to the roundoff error, but the results still suffer from inaccuracies due to the low number of sampling directions, because only a single orthogonal direction is used.
\item \textbf{Multi-Vector methods:} These are not expected to be as good as multi-axis methods for estimating gradients. However, they have the significant advantage that we can use any semi-regular polyhedra or semi-regular polytops with congruent vertices to obtain uniformly distributed vectors on the surfaces of unit spheres or unit hyperspheres.  
\end{enumerate}

\section{Conclusion and Future Work}

This is a preliminary study that shows the potential of the derivative estimation of Squire and Trapp for estimating the gradient in higher dimensions. The extension of these gradient estimation methods to higher dimensions is theoretically straightforward. Unfortunately, it is hard to compute the quality of the results, since ground truths for the gradients may not necessarily be reliable. We are still able to make some observations on the nature of these operations. We believe that there is strong merit in combining Squire and Trapp with the multi-axis extension of multi-sample method\cite{hart1989}. For example, we suspect that these methods can be used to obtain more accurate and fast solutions for optimization problems. Therefore, we believe that these methods could be better tested to solve optimization problems by demonstrating faster and more accurate results.


\bibliographystyle{unsrtnat}
\bibliography{references}

\end{document}